\documentclass{article}
\usepackage {amscd}
\usepackage {amsmath}
\usepackage {amssymb}
\usepackage [backend=biber,style=alphabetic,maxalphanames=4,sorting=nyt]{biblatex}
\usepackage {mathrsfs}
\usepackage {geometry}
\usepackage {setspace}
\usepackage {tikz}
\usepackage {tikz-cd}
\usepackage {float}
\usepackage {indentfirst}
\usepackage {enumitem}
\usepackage [amsmath,amsthm,thmmarks,hyperref,thref]{ntheorem}
\usepackage {hyperref}
\usepackage {mathabx}
\usepackage {diagbox}
\usepackage {graphicx}
\usepackage{longtable}
\usepackage{etoolbox}
\preto\tabular{\setcounter{magicrownumbers}{0}}
\preto\longtable{\setcounter{magicrownumbers}{0}}
\newcounter{magicrownumbers}

\graphicspath{ {./images/} }

\theoremstyle{plain}
\newtheorem{theorem}{Theorem}[section]
\newtheorem{proposition}[theorem]{Proposition}
\newtheorem{corollary}[theorem]{Corollary}
\newtheorem{lemma}[theorem]{Lemma}
\newtheorem{conjecture}[theorem]{Conjecture}

\theoremstyle{definition}
\newtheorem{definition}[theorem]{Definition}
\newtheorem{proposition-definition}[theorem]{Proposition-Definition}

\newtheorem{convention}[theorem]{Convention}

\theoremstyle{remark}

\newtheorem{construction}[theorem]{Construction}
\newtheorem{statement}[theorem]{}
\newtheorem*{remark}{Remark}

\newcolumntype{C}{>{$}c<{$}} 

\author {Yang He}
\title {On the strong Sarkisov program}

\begin {document}\large
\begin {spacing}{1.5}
\maketitle

\begin{abstract}
In this paper, we showed that strong Sarkisov program in dimension $d$ can be derived from termination of specific log flips in dimension $\leq d-1$. As a corollary, we prove that the strong Sarkisov program holds in dimension 4. Additionally, we prove the weak Sarkisov program with decreasing (augmented) Sarkisov degree. Finally, using the theory of syzygies of Mori fibre spaces, we obtain a directed diagram that encodes the information of the strong Sarkisov program.
\end{abstract}

\tableofcontents

\section{Introduction}

\subsection{History}
The study of the birational automorphism groups of the projective spaces $\mathbb{P}^{n}$ dates back to the 19th century, with significant contributions from the Italian school of algebraic geometry and notable figures like Emmy Noether. These groups, known as Cremona groups after Luigi Cremona, who initially introduced them, have been foundational in the field. One of the most famous theorems during this period is the classical theorem of Noether and Castelnuovo:

\begin{theorem}[Noether-Castelnuovo]
  The birational automorphism group of $\mathbb{P}^2$ is generated by the regular automorphism group and one quadratic transformation.
\end{theorem}

In the 1980s the minimal model program (MMP) was introduced, expanding the study of Cremona groups to encompass birational maps between Mori fibre spaces. Sarkisov, in his work \cite{Sarkisov89}, introduced the notion of \emph{elementary links} between Mori fibre spaces, and proposed a proof that every birational transformation between 3-fold Mori fibre spaces can be factorized into these elementary links. More precisely, he invented a procedure to explicitly untwist any birational map between polarized Mori fibre spaces by an elementary link, and introduced an invariant, the Sarkisov degree, that strictly decrease (with some exceptions) in this procedure. He stated that this procedure must terminate after finitely many steps. This procedure was later named after Sarkisov and called \emph{the Sarkisov program}. To distinguish with later developments, the termination of the procedure introduced by Sarkisov is referred to as the \emph{strong Sarkisov program}, c.f. \thref{construction of strong Sarkisov link} and \thref{conj:strong Sarkisov program}.

The first detailed proof of the strong Sarkisov program in dimension 3 was given by Corti:
\begin{theorem}[\cite{Cor}]
  The strong Sarkisov program holds in dimension 3.
\end{theorem}

Corti's proof was later summarized and simplified in \cite{BrunoMatsuki}, leading to the following conclusion:

\begin{theorem}[cf. \cite{BrunoMatsuki}]\label{theorem: BrunoMatsuki}
  Assume the following conditions hold:
  \begin{enumerate}
    \item the Log MMP conjecture holds in dimension $n$;
    \item the BAB conjecture holds in dimension $d$ for all $d \leq n$, that is, the family of $\epsilon$-klt $\mathbb{Q}$-factorial Fano varieties is bounded;
    \item the $S_{n}$(Local) property holds, that is, the local log canonical thresholds of $\mathbb{Q}$-factorial $n$-folds satisfy the ACC (ascending chain condition).
  \end{enumerate}
  Then the strong Sarkisov program holds in dimension $n$.
\end{theorem}

    The BAB conjecture was later proved by Birkar (cf. \cite{BAB}). The ACC for log canonical thresholds for usual pairs was proved in \cite{HMX}. The Log MMP conjecture, however, remains open in dimension $\geq 4$.

    Instead of looking into the construction of Sarkisov and tracing the Sarkisov degree, Hacon and McKernan used the theory of geography of log models to prove that the factorization into Sarkisov links always exists:
\begin{theorem}[\cite{HM}]
  Every birational map $f: X/S \dashrightarrow Y/T$ between Mori fibre spaces is a composition of Sarkisov links.
\end{theorem}

    However, the Sarkisov degree may not decrease for such factorization. We refer to this result as \emph{the weak Sarkisov program}.

Roughly speaking, the strong Sarkisov program is the combination of the minimal model program together with several other procedures. We can find the analog in the minimal model program as in the following table:

\[
\begin{matrix}
  \textrm{Sarkisov program} & \longleftrightarrow & \textrm{Minimal Model Program(MMP)}\\
  \textrm{Sarkisov links} & \longleftrightarrow & \textrm{divisorial contractions and flips}\\
  \textrm{Sarkisov degree} & \longleftrightarrow & \textrm{invariants in MMP}\\
  \textrm{weak Sarkisov program} & \longleftrightarrow & \textrm{existence of minimal model}\\
  \textrm{strong Sarkisov program} & \longleftrightarrow & \textrm{termination of MMP}
\end{matrix}
\]

\subsection{Main results}

The main goal of this article is to relax the general Log MMP assumption in \thref{theorem: BrunoMatsuki} and establish the strong Sarkisov program from termination of some specific types of Log MMP, i.e. to make precise the correspondence
$$
    \textrm{strong Sarkisov program}  \longleftrightarrow  \textrm{termination of MMP}.
$$

    To state the main results, we start by outlining several conjectures regarding the termination of flips. They are special cases of a more general conjecture about the termination of flips.

\begin{conjecture}[termination of canonical flips with big mobile boundary]\label{termination of canonical flips}
    Let $X$ be a $\mathbb{Q}$-factorial variety with terminal singularities and $\mathcal{H}$ be a mobile linear system of big divisors. Let $r$ be a positive rational number such that the pair $(X,r\mathcal{H})$ is canonical. Then any sequence of $(K_X + r\mathcal{H})$-flips terminates. More precisely, there is no infinite sequence of log flips
  $$
  (X,r\mathcal{H}) \dashrightarrow (X_{1},r\mathcal{H}_{X_{1}}) \dashrightarrow \cdots.
  $$
\end{conjecture}

\begin{conjecture}[termination of $\delta$-lc flips with big boundary]\label{termination of delta lc flips}
  Let $X$ be a $\mathbb{Q}$-factorial variety and $\Delta$ be a big boundary such that the pair $(X,\Delta)$ is $\delta$-lc. Then any sequence of $(K_X + \Delta)$-flips terminates. More precisely, there is no infinite sequence of log flips
  $$
  (X,\Delta) \dashrightarrow (X_{1},\Delta_{X_{1}}) \dashrightarrow \cdots.
  $$
\end{conjecture}

We also use the following conjecture of Shokurov which was recently proved by Birkar.

\begin{theorem}[Shokurov's $\epsilon-\delta$ conjecture for fibrations, cf. \cite{Birkar23}, Theorem 1.1]\label{Shokurov epsilon delta conjecture}
  Let $d$ be a positive integer and $\epsilon$ be a positive real number. Then there is a positive real number $\delta = \delta(d, \epsilon)$ satisfying the following:

  For any contraction $\pi : (X, B) \rightarrow Z$ such that
  \begin{itemize}[align=left]
    \item $\mathrm{dim}X - \mathrm{dim} Z = d$,
    \item $(X/Z, B)$ is $\epsilon$-lc,
    \item $K_{X}+B \sim_{\mathbb{Q},Z} 0$, and
    \item $-K_{X}$ is big over $Z$,
  \end{itemize}

    we can choose a divisor $M_{Z} \geq 0$ representing the moduli part, such that $(Z \ni z, B_{Z} +M_{Z})$ is $\delta$-lc.
\end{theorem}

The main theorem of this article is the following:

\begin{theorem}[also \thref{termination of strong Sarkisov link}]
  Assume one of the following conditions holds:
  \begin{enumerate}
    \item \thref{termination of canonical flips} for dimension $d$; or
    \item \thref{termination of delta lc flips} for $\delta = \delta(d-d',1-\epsilon)$ as in \thref{Shokurov epsilon delta conjecture} for all dimension $d' \leq d-1$ and some positive real number $\epsilon$.
  \end{enumerate}
  Then the strong Sarkisov program holds in dimension $d$.
\end{theorem}

    As a corollary, we prove that:

\begin{corollary}[also \thref{cor:strong Sarkisov program in dim 3 and 4}]
  The strong Sarkisov program holds in dimension 4.
\end{corollary}

    In the strong Sarkisov program, the Sarkisov degree may not decrease in certain exceptional cases. To fix this, we augment the Sarkisov degree, cf. \thref{def:augmented Sarkisov degree}. We have:

  \begin{proposition}[also \thref{prop:decrease Sarkisov degree}]
    The augmented Sarkisov degree decreases in every step of a strong Sarkisov program.
  \end{proposition}
  
  By making a specific choice in every step of the strong Sarkisov program, we are able to prove the weak Sarkisov program with decreasing augmented Sarkisov degree:
  
  \begin{theorem}[cf. \thref{theorem: weak Sarkisov program with decreasing Sarkisov degree} and \thref{cor:factorization with decreasing degree}]
    Let $\Phi: X/S \dashrightarrow Y/T$ be a birational map between Mori fibre spaces. Then there exists a terminate strong Sarkisov program for $\Phi$. In particular, there exists a factorization of $\Phi$ into Sarkisov links, such that the augmented Sarkisov degree decreases in every step.
  \end{theorem}

Using the theory of syzygies of Mori fibre spaces, we can construct a directed graph that encodes all the information of the strong Sarkisov program:
  \begin{theorem}[also \thref{thm:Sarkisov diagram}]
    Let $X'/S'$ be a Mori model of $Y$. Fix a very ample complete linear system $\mathcal{H}' = |m(-K_{X'} + f^{\prime *}A')|$ on $X'$, where $m$ is a sufficiently large integer and $A'$ is a very ample divisor on $S'$. Then for any Mori model $X/S$ of $Y$, there exists a directed graph $\mathfrak{G}(X/S,X'/S',\mathcal{H}')$ such that:
    \begin{enumerate}[label=(\arabic*)]
      \item The undirected graph of $\mathfrak{G}(X/S,X'/S',\mathcal{H}')$ is a finite subgraph of $CW_{1}(Y)$.
      \item The vertex corresponding to $X/S$ is the unique source, and the vertex corresponding to $X'/S'$ is the unique destination of $\mathfrak{G}(X/S,X'/S',\mathcal{H}')$.
      \item There is a 1-1 correspondence
      \begin{align*}
      \{\text{Paths of } \mathfrak{G} \text{ from the vertex corresponding to } X/S \text{ to the vertex corresponding to } X'/S'\} \\
       \longleftrightarrow \{ \text{Sequences of untwisting in the strong Sarkisov program of } X/S \dashrightarrow X'/S'\}
      \end{align*}
      and an injection
      \begin{align*}
      \{ \text{Steps of intermediate untwisting in the strong Sarkisov program of } X/S \dashrightarrow X'/S' \} \\
       \hookrightarrow \{\text{Edges of }\mathfrak{G}(X/S,X'/S',\mathcal{H}') \}.
      \end{align*}
    \end{enumerate}
  \end{theorem}

  \section*{Acknowledgement}

  The author wants to thank professor Vyacheslav V. Shokurov for proposing the question as well as several useful comments to this article.

\section{Preliminaries}

\subsection{Pairs with linear systems}
\begin{definition}[$\mathbb{R}$-linear system]
  Let $X$ be a $\mathbb{Q}$-factorial variety and $\mathcal{H}$ be an integral linear system on $X$. Let $c \geq 0$ be a real number. Define
  $$
  c\mathcal{H}_{\mathbb{R}} = \{ \sum\limits_{i}b_{i}H_{i}\mid H_{i} \in \mathcal{H},\ b_{i} \geq 0,\ \sum\limits_{i}b_{i}=c \}.
  $$
  We say that $c\mathcal{H}_{\mathbb{R}}$ is an \emph{$\mathbb{R}$-linear system}. Similarly, we define
  $$
  c\mathcal{H} = \{ cH \mid H \in \mathcal{H}\}.
  $$
  We say that $c\mathcal{H}$ is an \emph{integral linear system scaled by $c$}. In particular, if $c=1$, it becomes a linear system in usual sense.
\end{definition}

\begin{proposition}
  Let $F$ be the fixed part of the linear system $\mathcal{H}$. Then for every $D \in c\mathcal{H}_{\mathbb{R}}$ we have $D \geq cF$, and $cF$ is the largest divisor satisfying this property, that is, if for every $D \in c\mathcal{H}_{\mathbb{R}}$ we have $D \geq F'$ for some $\mathbb{R}$-divisor $F'$, then $cF \geq F'$.
\end{proposition}

\begin{proof}
  It is obvious by definition that for every $D \in c\mathcal{H}_{\mathbb{R}}$ we have $D \geq cF$. It suffices to show that $cF$ is the largest one. Indeed, if $F'$ is another such divisor satisfying this property, then for any divisor $H \in \mathcal{H}$, we have $cH \geq F'$. Since $F$ is the fixed part of $\mathcal{H}$, we conclude that $cF \geq F'$.
\end{proof}

\begin{statement}[Fixed and mobile part]
  We say that $cF$ is the \emph{($\mathbb{R}$-)fixed part} of the $\mathbb{R}$-linear system $c\mathcal{H}_{\mathbb{R}}$. Write $\mathcal{H} = \mathcal{M} + F$. Then we say that the $\mathbb{R}$-linear system $c\mathcal{M}_{\mathbb{R}}$ is the \emph{($\mathbb{R}$-)mobile part} of $c\mathcal{H}_{\mathbb{R}}$.
\end{statement}

\begin{definition}[Linear systems and log discrepancies]
  Let $\mathcal{H}$ be an $\mathbb{R}$-linear system or an integral linear system scaled by some real number on a $\mathbb{Q}$-factorial variety $X$. Let $x \in X$ be a point. Then we can define the local log discrepancy
  $$
  a(x \in X, \mathcal{H}) = \sup\limits_{B \in \mathcal{H}} a(x \in X, B).
  $$
  For an exceptional divisor $E$ over $X$, we define the log discrepancy
  $$
  a(E, X, \mathcal{H}) = \sup\limits_{B \in \mathcal{H}} a(E, X, B).
  $$
  Globally, we define
  $$
  a(X, \mathcal{H}) = \inf\limits_{x \in X} a(x \in X, \mathcal{H}) = \inf\limits_{E\text{ exceptional}} a(E, X, \mathcal{H}).
  $$
  We say that the pair $(X,\mathcal{H})$ is
  $$
  \begin{cases}
    \text{log canonical (lc)}, & \mbox{if } a(X, \mathcal{H}) \geq 0, \\
    \text{canonical}, & \mbox{if } a(X, \mathcal{H}) \geq 1, \\
    \text{terminal}, & \mbox{if } a(X, \mathcal{H}) >1.
  \end{cases}
  $$
  We say that the pair $(X,\mathcal{H})$ is \emph{Kawamata log terminal (klt)}, if $a(X, \mathcal{H}) > 0$, and the $\mathbb{R}$-fixed part of $\mathcal{H}$ has coefficients $ < 1$. For some positive real number $\delta$, we say that the pair $(X,\mathcal{H})$ is $\delta$-lc (resp. $\delta$-klt), if $a(X, \mathcal{H}) \geq \delta$ (resp. $> \delta$), and the $\mathbb{R}$-fixed part of $\mathcal{H}$ has coefficients $ \leq 1 - \delta$ (resp. $ < 1-\delta$).
\end{definition}

\begin{proposition}
  Let $X$ be a $\mathbb{Q}$-factorial variety and $\mathcal{H}$ a linear system on $X$. Then we have
  $$
  a(X,c\mathcal{H}) = a(X,c\mathcal{H}_{\mathbb{R}})
  $$
  for any $0 \leq c \leq 1$.
\end{proposition}

\begin{proof}
    Let $f:Y \rightarrow X$ be a log resolution such that $f^{*}\mathcal{H} = \mathcal{M} + F$, where the fixed part $F$ has snc support and $\mathcal{M}$ is base point free. Then
  $$
  a(X,c\mathcal{H}) = \min\{a(Y,cF) , \min\limits_{E_{i}\ f-\text{exceptional}}\{ a(E_{i},X,c\mathcal{H}) \}\} = a(X,c\mathcal{H}_{\mathbb{R}}).
  $$
\end{proof}

The above proposition shows that the notion of $\mathbb{R}$-linear systems doesn't introduce ambiguity for log discrepancies. By abusing notations, we may remove the subscript and also denote the $\mathbb{R}$-linear system $\mathcal{H}_{\mathbb{R}}$ by $\mathcal{H}$.

For a pair involving $\mathbb{R}$-linear systems, we often use the following approximation result to reduce the discussion to the discussion for usual pairs:

\begin{proposition}\label{prop:approach ld by divisors}
  Let $(X,\mathcal{H})$ be a pair where $\mathcal{H}$ is an $\mathbb{R}$-linear system. Then for any $\epsilon > 0$, there exists a divisor $D \in \mathcal{H}$ such that $a(X,\mathcal{H}) - \epsilon < a(X,D) \leq a(X,\mathcal{H})$, and the mobile part of $\mathcal{H}$ has coefficients $< \epsilon$. Moreover, if $\mathcal{H}$ is mobile, then we can let $a(X,D) = a(X,\mathcal{H})$.
\end{proposition}

\begin{proof}
  Let $f:Y \rightarrow X$ be a log resolution such that $f^{*}\mathcal{H} = \mathcal{M} + F$, where the fixed part $F$ has snc support and the $\mathbb{R}$-linear system $\mathcal{M}$ is base point free. Write $F = \sum f_{i}F_{i}$. Then
  $$
  a(X,\mathcal{H}) = \min\{a(Y,F) , \min\limits_{F_{i}\ f-\text{exceptional}}\{ 1 - f_{i}\}\}.
  $$
  Choose $M \in \mathcal{M}$ such that
  \begin{enumerate}
    \item $M$ is snc;
    \item $M$ contains no components of $F$ and no $f$-exceptional components;
    \item $M$ doesn't pass through the center of exceptional divisors of log discrepancy 
        $$
        < totallogdiscrepancy(Y,F) + 1;
        $$
    \item the coefficients of $M$ is $< \frac{\epsilon}{\mathrm{dim}X}$.
  \end{enumerate}
  Let $D \in \mathcal{H}$ be the $\mathbb{R}$-divisor such that $f^{*}D = M + F$. Then the divisor $D$ satisfies our requirement.
\end{proof}

The following proposition shows the benefits of using ($\mathbb{R}$-)linear systems instead of usual divisors:

\begin{proposition}\label{prop:singularity after MMP}
  Let $(X,\mathcal{H})$ be a terminal (resp. canonical) pair, where $\mathcal{H}$ is a mobile $\mathbb{R}$-linear system. Let $g:X \dashrightarrow X'$ be a $(K_{X}+\mathcal{H})$-negative divisorial contraction or a log flip. Then the pair $(X',\mathcal{H}_{X'})$ is still terminal (resp. canonical).
\end{proposition}

\begin{proof}
  Since $\mathcal{H}$ is mobile, we can choose a general element $D \in \mathcal{H}$ such that $g$ doesn't contract prime components of $D$, and
  $$
  a(X,\mathcal{H})-\epsilon < a(X,D) \leq a(X,\mathcal{H})
  $$
  for some positive real number $\epsilon \ll 1$. Then for any divisor $E$ on $X$ contracted by $g$, we have
  $$
  a(E;X',D_{X'}) > 1.
  $$
  Since $g$ is a $(K_{X}+\mathcal{H})$-negative divisorial contraction or log flip, we have
  $$
  a(X',D_{X'}) \geq \min \{a(X,D),a(E;X',D_{X'}) \} > \min \{a(X,\mathcal{H})-\epsilon , 1\}.
  $$
  Hence the pair $(X',\mathcal{H}_{X'})$ is still terminal (resp. canonical).
\end{proof}

\subsection{Preliminaries about the strong Sarkisov program}
\begin{convention}
  In this article, unless specified, every contraction $X \rightarrow S$ is a contraction between \emph{projective varieties}. A Mori fibre space is $\mathbb{Q}$-factorial and has terminal singularities.
\end{convention}

\begin{definition}[Sarkisov degree]\label{Sarkisov degree}
  Let
  $$
  \begin{tikzcd}
    X \arrow[dashed,r,"\Phi"]\arrow[d,"f"] & X' \arrow[d,"f'"] \\
    S & S'
  \end{tikzcd}
  $$
  be a birational map between Mori fibre spaces. Fix a very ample complete linear system $\mathcal{H}' = |m(-K_{X'} + f^{\prime *}A')|$ on $X'$, where $m$ is a sufficiently divisible integer and $A'$ is a very ample divisor on $S'$, and denote $\mathcal{H}$ the total transform of $\mathcal{H}'$ on $X$, that is we take a common resolution
  $$
  \begin{tikzcd}
    & W \arrow[ld,"p"] \arrow[rd,"q"] & \\
    X \arrow[dashed,rr]& & X'
  \end{tikzcd}
  $$
  
    and let $\mathcal{H} = p_{*}(q^{*}\mathcal{H}')$. For sufficiently divisible $m$, divisors in $\mathcal{H}'$ are Cartier and $\mathcal{H}$ is an integral linear system. Moreover, the linear system $\mathcal{H}$ doesn't depend on the choice of the common resolution $W$. Then the \emph{Sarkisov degree} $\mathrm{deg}(\Phi,\mathcal{H}')$ is the triple $(\mu,c,e)$ defined as follows:

  \begin{enumerate}[label=(\arabic*)]
    \item The invariant $\mu \in \mathbb{Q}_{+}$ is the unique number such that
     $$
     \mathcal{H} + \mu K_X \sim_{\mathbb{R},S} 0
     $$
     The number exists since $X/S$ is a Mori fibre space, and it's positive since $\mathcal{H}$ is ample over $S$. Notice that the reciprocity $\frac{1}{\mu}$ is the \emph{pseudo-effective threshold} of $X/S$ with respect to $\mathcal{H}$.
    \item The invariant $c \in \mathbb{Q}_{+}$ is the canonical threshold of $X$ with respect to $\mathcal{H}$. That is,
    $$
    c = \mathrm{sup}\{t \in \mathbb{R}|\text{the pair }(X,t\mathcal{H})\text{ is canonical}\}.
    $$
    \item The invariant $e \in \mathbb{Z}_{\geq 0}$ is the number of crepant exceptional divisors over the pair $(X, c\mathcal{H})$.
  \end{enumerate}

    The \emph{order of the Sarkisov degree} is defined as follows: $(\mu_1,c_1,e_1) > (\mu_2,c_2,e_2)$ if either:

  \begin{enumerate}[label=(\alph*)]
    \item $\mu_1 > \mu_2$, or
    \item $\mu_1 = \mu_2$ and $c_1 < c_2$, or
    \item $\mu_1 = \mu_2$, $c_1 = c_1$ and $e_1 > e_2$.
  \end{enumerate}

\end{definition}

\begin{remark}
  In general, $\mathcal{H}$ may not be a complete linear system.
\end{remark}

  \begin{definition}[square map and square isomorphism]\label{square map}
    A birational map $\Phi: X/S \dashrightarrow X'/S'$ between Mori fibre spaces is called \emph{square} if there exists a \emph{birational} map $\phi: S \dashrightarrow S'$ such that the following diagram is commutative
    $$
    \begin{tikzcd}
      X \arrow[r,dashed,"\Phi"]\arrow[d] & X' \arrow[d] \\
      S \arrow[r,dashed,"\phi"]      &  S'
    \end{tikzcd}
    $$
    and $\Phi$ is an isomorphism over the generic point of $S$. We say that $\Phi$ is a \emph{square isomorphism}, if $\Phi$ is an isomorphism, and the induced map $\phi$ is also an isomorphism.
  \end{definition}

\begin{theorem}[Noether-Fano-Iskovskikh inequality, cf. \cite{Isk04}, Theorem 3.3]\label{Noether-Fano-Iskovskikh inequality}
  Let
  $$
  \begin{tikzcd}
    X \arrow[dashed,r,"\Phi"] \arrow[d,"f"] & X' \arrow[d,"f'"] \\
    S & S'
  \end{tikzcd}
  $$
  be a birational map between Mori fibre spaces. Let $\mathcal{H}'$ be a very ample linear system on $X'$ such that $K_{X'} + \frac{1}{\mu'}\mathcal{H}' \sim_{\mathbb{R}} f^{\prime *} A$ for some very ample linear system $A$ on $S'$, and $\mathcal{H}$ its total transformation on $X$. Let $(\mu,c,e)$ be the Sarkisov degree of $\Phi$ with respect to the linear system $\mathcal{H}'$, and $(\mu',c',e')$ the Sarkisov degree of the identity map of $X'/S'$ with respect to the linear system $\mathcal{H}'$. Then
  \begin{enumerate}
    \item We have the inequality $\mu \geq \mu'$, and equality holds only if $\Phi$ induces a rational map $\phi:S \dashrightarrow S'$, such that the following diagram is commutative:
        $$
        \begin{tikzcd}
        X \arrow[dashed,r,"\Phi"] \arrow[d,"f"] & X' \arrow[d,"f'"] \\
        S \arrow[dashed,r,"\phi"] & S'
        \end{tikzcd}
        $$
    \item Assume further that the pair $(X,\mu \mathcal{H})$ is canonical and $K_X+\frac{1}{\mu}\mathcal{H}$ is (absolutely) nef, then $\Phi$ is a square isomorphism. In particular, $\mu = \mu'$.
  \end{enumerate}

\end{theorem}

\begin{corollary}\label{not square isomorphism}
  If $\Phi$ is not a square isomorphism, then either:

  \begin{enumerate}[label=(\arabic*)]
    \item The pair $(X,\mu \mathcal{H})$ is not canonical, or
    \item The pair $(X,\mu \mathcal{H})$ is canonical, but $K_X+\frac{1}{\mu}\mathcal{H}$ is not (absolutely) nef.
  \end{enumerate}
\end{corollary}

\begin{proof}
  This is straightforward by \thref{Noether-Fano-Iskovskikh inequality}.
\end{proof}

\begin{definition}\label{def:maximal extraction and maximal center}
   Let settings be as in \thref{Sarkisov degree} and assume that $c < \frac{1}{\mu}$.

  \begin{enumerate}[label=(\roman*)]
    \item A \emph{maximal extraction} is an extremal divisorial extraction $f: \widetilde{X} \rightarrow X$ of a divisor $E$ of log discrepancy $a(E,X,c\mathcal{H}) = 1$ such that $\widetilde{X}$ is $\mathbb{Q}$-factorial with terminal singularities.
    \item An exceptional divisor $E$ over $X$ is called a \emph{maximal singularity} of $\mathcal{H}$ if $E$ is the exceptional divisor of a maximal extraction. The center $\mathrm{Center}_X(E)$ is called the \emph{maximal center}.
  \end{enumerate}

\end{definition}

\begin{proposition}[cf. \cite{Cor}, Proposition-Definition 2.10]\label{prop:existence of maximal extraction}
  Let notations be as in \thref{Sarkisov degree}. Assume that $c < \frac{1}{\mu}$. Then $X$ admits a maximal singularity $E$, i.e. there exists a maximal extraction of some exceptional divisor $E$.
\end{proposition}

\begin{definition}[Sarkisov links]\label{Sarkisov link}
  A Sarkisov link between Mori fibre spaces $\phi: X \rightarrow S$ and $\psi: Y \rightarrow T$ is a diagram of birational maps with one of the following types:
  
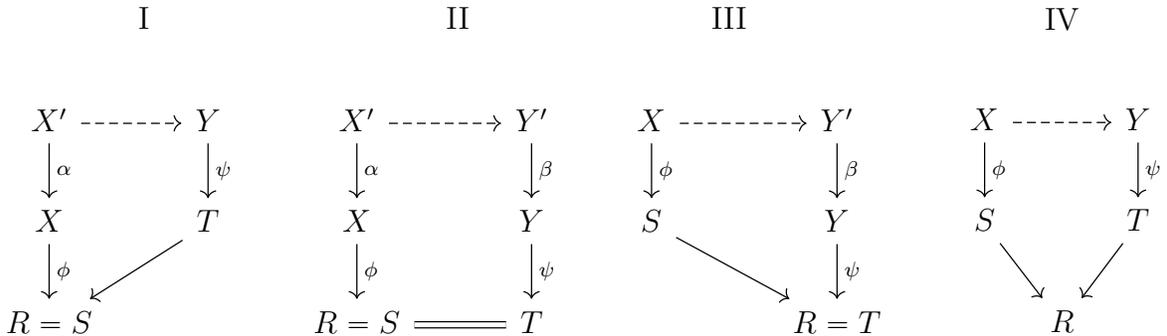
\begin{equation*}
    \begin{tikzcd}[column sep=tiny]
        &\mathrm{I}&\\
        X' \arrow[rr,dashed] \arrow[d,"\alpha"] && Y \arrow[d,"\psi"] \\
        X \arrow[d,"\phi"] && T \arrow[lld]\\
        R=S
    \end{tikzcd}
    \qquad
    \begin{tikzcd}[column sep=tiny]
         &\mathrm{II}&\\
        X' \arrow[rr,dashed] \arrow[d,"\alpha"] && Y' \arrow[d,"\beta"] \\
        X \arrow[d,"\phi"] && Y \arrow[d,"\psi"]\\
        R=S \arrow[rr,equal] && T
    \end{tikzcd}
    \qquad
    \begin{tikzcd}[column sep=tiny]
        &\mathrm{III}&\\
        X \arrow[rr,dashed] \arrow[d,"\phi"] && Y' \arrow[d,"\beta"] \\
        S \arrow[rrd] && Y \arrow[d,"\psi"]\\
         &&R=T
    \end{tikzcd}
    \qquad
    \begin{tikzcd}[column sep=tiny]
        &\mathrm{IV}&\\
        X \arrow[rr,dashed] \arrow[d,"\phi"] && Y \arrow[d,"\psi"] \\
        S \arrow[rd] && T \arrow[ld]\\
        &R&
    \end{tikzcd}
  \end{equation*}
  The horizontal dashed arrows are small transformations of terminal varieties over $R$. All non-dashed arrows are extremal contractions, where $\alpha$ and $\beta$ are divisorial contractions.
\end{definition}

\begin{remark}
    It was pointed out by Shokurov that Sarkisov links can also be defined using the notion of central models, c.f. \cite{Myself_Syzygy}, Definition 5.2. 
\end{remark}

\subsection{Adjunction}

We need the following lemma for taking adjunction on a general fibre:

\begin{lemma}\label{adjunction to fibre}
  Let $f:X \rightarrow Z$ be a proper morphism between quasi-projective varieties. Suppose that $X$ is $\mathbb{Q}$-factorial with terminal singularities. Let $(X/Z,B)$ be a klt pair. Then for a general fibre $F$ of $f$, there exists a boundary $B_{F}$ such that $(F,B_{F})$ is klt and $(K_{X}+B)|_{F} = K_{F}+B_{F}$, where $B_F = B|_{F}$.
\end{lemma}

\begin{proof}
  We show that we can find a klt pair after taking a general hyperplane section on $Z$. Indeed, since $(X/Z,B)$ is klt, there are only finitely many exceptional divisors $E$ on $X$ with log discrepancy $0 < a(E;X,B) \leq 1$. Moreover, there are only finitely many centers on $X$ of multiplicity $>1$. Hence we can find a sufficiently general hyperplane section $H$ on $Z$ such that $f^{*}H$ doesn't contain these centers and is reduced and irreducible. Then the pair $(X/Z,B+f^{*}H)$ is plt. By adjunction, the pair $(f^{*}H,B_{f^{*}H})$ is klt. Moreover, since $X$ has terminal singularities, it is smooth in codimension 2. Hence for a general hyperplane $H$ we know that $f^{*}H$ is also smooth in codimension 2, so there is no difference part and $(K_{X}+B)|_{f^{*}H} = K_{f^{*}H} + B_{f^{*}H}$ where $B_{f^{*}H} = B|_{f^{*}H}$. Repeat this procedure and we obtain the result.
\end{proof}

\subsection{ACC and finiteness}

We state some important results about ACC and finiteness in the minimal model program, but in the language of pairs with linear systems.

\begin{theorem}[ACC for log canonical thresholds, cf. \cite{HMX}, Thm 1.1]\label{ACC for log canonical thresholds}
  Fix a positive integer $d$. Then the set
  $$
  LCT(d) = \{ \mathrm{lct}(X;\mathcal{H}) | \mathrm{dim}(X)=d,\ X \text{ is }\mathbb{Q}\text{-factorial lc},\ \mathcal{H} \text{ a mobile (integral) linear system} \}
  $$
   satisfies ACC.
\end{theorem}

\begin{proof}
  Take a resolution $f: Y \rightarrow X$ such that $f^{*}\mathcal{H} = \mathcal{M} + F$ where $\mathcal{M}$ is a base-point free linear system and the fixed part $F$ has simple normal crossing. Suppose
  $$
  f^{*}K_{X} = K_{Y} + \sum\limits_{i} a_{i}E_{i}; \qquad F = \sum\limits_{i}b_{i}E_{i}
  $$
  Let $M \in \mathcal{M}$ be a divisor such that the $D=M+F$ has simple normal crossing, and the components of $M$ has coefficient 1. Then for $0 \leq t \leq 1$, the pair $(X,t\mathcal{H})$ is log canonical if and only if the pair $(X,tD_{X})$ is log canonical. Hence we reduce to ACC for log canonical thresholds of normal pairs, i.e. \cite{HMX}, Thm 1.1.
\end{proof}

\begin{corollary}[ACC for Fano index, cf.  \cite{HMX}, Corollary 1.10]\label{ACC for Fano index}
  Let $I$ be a DCC set and $d$ a positive integer. Let $\mathfrak{D}$ be the set of projective log canonical pairs $(X,\Delta)$ of dimension $d$ such that $-(K_{X}+\Delta)$ is ample and the coefficients of $\Delta$ belong to $I$. The \emph{Fano Index} of $(X,\Delta)$ is the largest real number $r$ such that we can write $-(K_{X}+\Delta) \sim_{\mathbb{R}} rH$ for some ample Cartier divisor $H$. We define the set of Fano indexes of Fano varieties of dimension $d$ to be
  $$
  R = R_{d}(I) = \{ r \in \mathbb{R} \mid r\text{ is the Fano index of }(X,\Delta) \in \mathfrak{D}  \}
  $$
  Then $R$ satisfies ACC.
\end{corollary}

\begin{remark}
    The BAB conjecture assumption in \thref{theorem: BrunoMatsuki} is only used to derive \thref{ACC for Fano index}. Hence the proof actually doesn't rely on the BAB conjecture.
\end{remark}

\begin{proposition}[Lower semi-continuity of canonical threshold]\label{lower semi-continuity of canonical threshold}
  Let $X$ be a variety with canonical singularities, $\mathcal{H}$ a mobile linear system on $X$. Then the function $x \mapsto c(x,X,\mathcal{H})$ is lower semi-continuous on $X$ in Zariski topology. In particular, there are only finitely many possible values of $c(x,X,\mathcal{H})$.
\end{proposition}

\begin{proof}
  Take a resolution $f: Y \rightarrow X$ such that $f^{*}\mathcal{H} = \mathcal{M} + F$ where $\mathcal{M}$ is a base-point free linear system and the fixed part $F$ has simple normal crossing. We have $\mathrm{Supp}F \subseteq \mathrm{exc}(f)$. Suppose
  $$
  f^{*}K_{X} = K_{Y} + \sum\limits_{i} a_{i}E_{i}; \qquad F = \sum\limits_{i}b_{i}E_{i}
  $$
  where $E_i$ are prime exceptional divisors, $a_{i} \leq 0$ and $b_{i} \geq 0$. Then
  $$
  a(x \in X,c\mathcal{H}) = \min\{2,\min_{x \in f(E_{i})} \{1-a_{i}-cb_{i}\}, \min_{x \in f(E_{j}\cap E_{k})} \{2-a_{j}-a_{k}-cb_{j}-cb_{k}\} \}.
  $$
  Hence
  $$
  c(x \in X, \mathcal{H}) = \min \{\min_{\substack{x \in f(E_{i}) \\ b_{i} > 0}} \{ \frac{-a_{i}}{b_{i}} \} , \min_{\substack{x \in f(E_{j}\cap E_{k}) \\ b_{i} + b_{j} > 0}}\{\frac{1-a_{j}-a_{k}}{b_{j}+b_{k}}\} \}.
  $$
  Hence it is lower semi-continuous.
\end{proof}

\begin{remark}
The same holds for any $a$-log canonical threshold.
\end{remark}

\begin{lemma}\label{finiteness of divisors with negative discrepancy}
  Let $(X,\alpha\mathcal{H})$ be a pair with only plt singularities, where $\mathcal{H}$ is a mobile linear system. Then there are only finitely many divisors $E$ over $X$ with $a(E,X,\alpha\mathcal{H}) < 1$.
\end{lemma}

\begin{proof}
  Take a resolution $f: Y \rightarrow X$ such that $f^{*}\mathcal{H} = \mathcal{M} + F$ where $\mathcal{M}$ is a base-point free linear system and the fixed part $F$ has simple normal crossing. Since $\mathcal{H}$ is mobile, we have $\mathrm{Supp}F \subseteq \mathrm{exc}(f)$. Suppose that
  $$
  f^{*}K_{X} = K_{Y} + \sum\limits_{i} a_{i}E_{i} \qquad \text{and} \qquad F = \sum\limits_{i}b_{i}E_{i}.
  $$
  Then $(X,\alpha\mathcal{H})$ being plt implies $(Y,\sum\limits_{i} (a_{i} + \alpha b_{i})E_{i})$ is klt. Hence the number of divisors with log discrepancy $<1$ is finite.
\end{proof}

\section{The strong Sarkisov program}\label{strong Sarkisov program}

\subsection{Construction of the untwisting}

We begin the strong Sarkisov program as follows:

\begin{construction}[2-ray game, cf. \cite{Cor}, Theorem 5.4 and \cite{Matsuki} \S 13.1]\label{construction of strong Sarkisov link}
  Let
  $$
  \begin{tikzcd}
    X \arrow[dashed,r,"\Phi"]\arrow[d,"f"] & X' \arrow[d,"f'"] \\
    S & S'
  \end{tikzcd}
  $$
  be a birational map between Mori fibre spaces. We define the Sarkisov degree $\mathrm{deg}(\Phi,\mathcal{H}') = (\mu,c,e)$ as in \thref{Sarkisov degree}. If $\Phi$ is a square isomorphism then we do nothing. Otherwise by \thref{not square isomorphism} we have:

  \begin{enumerate}[label = \textbf{Case \arabic*:}, align=left]

  \item $K_X + \frac{1}{\mu}\mathcal{H}$ is not canonical. In this case the canonical threshold $c < \frac{1}{\mu}$, and there exists a maximal extraction of $(X,c\mathcal{H})$ by \thref{prop:existence of maximal extraction}. Let $f: \widetilde{X} \rightarrow X/S$ be a maximal extraction and denote $(\widetilde{X},c\mathcal{H}_{\widetilde{X}})$ the crepant pullback of the pair $(X,c\mathcal{H})$. Since the extraction is maximal, $\mathcal{H}_{\widetilde{X}}$ is just the strict transformation of $\mathcal{H}$. Then $\rho(\widetilde{X}/S)=2$ and hence $\overline{NE(\widetilde{X}/S)}$ is generated by 2 extremal rays. One extremal ray $R$ is $(K_{\tilde{X}}+c\mathcal{H}_{\widetilde{X}})$-trivial and leads to $\widetilde{X} \rightarrow X/S$ by construction. Since $c < \frac{1}{\mu}$, the other extremal ray $R'$ is $(K_{\widetilde{X}}+c\mathcal{H}_{\widetilde{X}})$-negative. We run the Log MMP on the pair $(\widetilde{X}/S,c\mathcal{H}_{\widetilde{X}})$ to contract the extremal ray $R'$. There are 3 possible situations:

\begin{enumerate}[label=(\roman*), align=left]
  \item A log Mori fibre space;
  \item A log divisorial contraction, in this case we get a log Mori fibre space after this divisorial contraction;
  \item A log flip, in this case the new cone $\overline{NE}$ of effective curves will be generated by a positive extremal ray since we take a log flip, and a negative extremal ray since we start from a non-pseudo-effective pair. We then continue to contract the negative extremal ray and again we get one of the 3 situations. Since $\widetilde{X}/S$ has Fano type, a sequence of log flips terminates after finitely many steps, so after getting into situation (iii) finitely many times, we will finally get into situation (i) or (ii).
\end{enumerate}

In all situations we will get a log Mori fibre space after finitely many steps.

Since the resulting pairs are obtained from a $(K_{\widetilde{X}}+c\mathcal{H}_{\widetilde{X}})$-MMP over $S$, they have canonical singularities. Furthermore, since $\widetilde{X}$ is $\mathbb{Q}$-factorial with terminal singularities, every crepant exceptional divisor of the pair $(\widetilde{X}, c\mathcal{H}_{\widetilde{X}})$ must center at the base locus of $\mathcal{H}_{\widetilde{X}}$. Hence the new Mori fibre space $X''/S''$ has terminal singularities. In this case, the Sarkisov link is of type I or type II.

  \item The pair $(X,\frac{1}{\mu}\mathcal{H})$ is canonical but $K_X + \frac{1}{\mu}\mathcal{H}$ is not (absolutely) nef. Then $c \geq \frac{1}{\mu}$. Let $R$ be the extremal ray contracted in $X \longrightarrow S$, then $R$ is $(K_X + \frac{1}{\mu}\mathcal{H})$-trivial and $K_{X}$-negative. We claim that we can find a $(K_X + \frac{1}{\mu}\mathcal{H})$-negative extremal ray $R'$ such that $R$ and $R'$ span an extremal face. Indeed, since $\mathcal{H}$ is big, we can apply the cone theorem and know that there are only finitely many $(K_X + \frac{1}{\mu}\mathcal{H})$-negative and $(K_X + \frac{1}{\mu}\mathcal{H})$-trivial extremal rays. Let $L$ be a supporting hyperplane of $R$, that is, $L \cdot R = 0$ and $L \cdot v > 0$ for $v \in \overline{NE}(X) - R$. Let $R'$ be the $(K_X + \frac{1}{\mu}\mathcal{H})$-negative extremal ray such that the number $a = \frac{-(K_X + \frac{1}{\mu}\mathcal{H}) \cdot v'}{L \cdot v'} > 0$ is the largest. If we take a sufficiently general $L$, such extremal ray will be unique. Then the hyperplane defined by $(K_X + \frac{1}{\mu}\mathcal{H}+aL)\cdot v = 0$ is the supporting hyperplane of the face spanned by $R$ and $R'$.

      The extremal contraction of rays $R$ and $R'$ gives the following diagram:

      $$
      
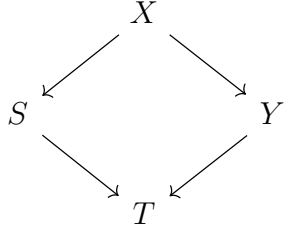
\begin{tikzcd}
         & X \arrow{ld} \arrow{rd}  & \\
         S \arrow{rd} & &    Y \arrow{ld} \\
         & T &
      \end{tikzcd}
      $$

      Similar to Case 1, we run the $(K_X + \frac{1}{\mu}\mathcal{H})$-MMP on $X/T$ such that the first extremal contraction is $X \rightarrow Y$. Hence we obtain the desired Sarkisov link. In this case, the Sarkisov link is of type III or type IV.
   \end{enumerate}

   The above construction is called a \emph{2-ray game}. We say that the Sarkisov link $X/S \dashrightarrow X''/S''$ constructed above is an \emph{untwisting} of the birational map $X/S \dashrightarrow X'/S'$.

   Finally, we replace $X/S$ by $X''/S''$ and repeat the construction above. We will obtain a sequence of Sarkisov links:
   $$
   X/S \dashrightarrow X_{1}/S_{1} \dashrightarrow X_{2}/S_{2} \dashrightarrow \cdots
   $$
\end{construction}

\begin{remark}
The constructions of the 2-ray game for Case 1 are slightly different between \cite{Cor} and \cite{Matsuki}. However, these constructions give exactly the same result, since the extremal contraction is unique at each step.
\end{remark}

\subsection{Termination of the sequence of untwisting}

The remaining part of the strong Sarkisov program is the termination of the above construction:

\begin{conjecture}[Strong Sarkisov program]\label{conj:strong Sarkisov program}
  The sequence of Sarkisov links constructed in \thref{construction of strong Sarkisov link} terminates. In other words, the birational map $X_{i}/S_{i} \dashrightarrow X'/S'$ becomes a square isomorphism for some $i$.
\end{conjecture}

\begin{theorem}\label{termination of strong Sarkisov link}
  Assume that one of the following conditions holds:
  \begin{enumerate}
    \item \thref{termination of canonical flips} for dimension $d$; or
    \item \thref{termination of delta lc flips} for $\delta = \delta(d-d',1-\epsilon)$ as in \thref{Shokurov epsilon delta conjecture} for all dimension $d' \leq d-1$ and some positive real number $\epsilon$.
  \end{enumerate}

  Then the above construction in dimension $d$ terminates in finitely many steps.

\end{theorem}

\begin{proof}
  The proof basically follows the same framework as in \cite{Cor} and \cite{BrunoMatsuki} for dimension 3, with some additional cases to discuss in higher dimension. We always denote the new Mori fibre space obtained from \thref{construction of strong Sarkisov link} by $X'' \rightarrow S''$.

  Suppose first that we are in Case 2 of \thref{construction of strong Sarkisov link} and obtain a Sarkisov link of type IV. Then either
  \begin{enumerate}[label=2-(IV\alph*):,align=left]
    \item $K_{X''} + \frac{1}{\mu}\mathcal{H}''$ is negative over $S''$, or
    \item $K_{X''} + \frac{1}{\mu}\mathcal{H}''$ is nef over $T$, and its Iitaka fibration over $T$ is exactly $X'' \rightarrow S''$.
  \end{enumerate}
  Similarly, for a Sarkisov link of type III, either

    \begin{enumerate}[label=2-(III\alph*):,align=left]
    \item $K_{X''} + \frac{1}{\mu}\mathcal{H}''$ is negative over $S''$, or
    \item $K_{X''} + \frac{1}{\mu}\mathcal{H}''$ is nef over $T$, and its Iitaka fibration over $T$ is exactly $X'' \rightarrow S''$.
  \end{enumerate}

  Indeed, suppose $K_{X''} + \frac{1}{\mu}\mathcal{H}''$ is nef over $T$. Notice that the $(K_{X}+\frac{1}{\mu}\mathcal{H})$-MMP is also a $(K_{X}+(\frac{1}{\mu}-\epsilon)\mathcal{H})$-MMP for $0 < \epsilon \ll 1$. Since $K_{X}+(\frac{1}{\mu}-\epsilon)\mathcal{H}$ is anti-ample over $T$ for $0 < \epsilon \ll 1$, the fibration $X'' \rightarrow S''$ must be $(K_{X}+(\frac{1}{\mu}-\epsilon)\mathcal{H})$-negative for $0 < \epsilon \ll 1$. Hence it is $(K_{X''} + \frac{1}{\mu}\mathcal{H}'')$-trivial.

  Now we begin to discuss all the possible cases.

  \begin{enumerate}[align=left]
   \item \textbf{Case 2-(IIIa) and 2-(IVa)}. Suppose we are in the former cases 2-(IVa) or 2-(IIIa). Then by the definition of Sarkisov degree, the pseudo-effective threshold $\frac{1}{\mu}$ strictly increase in this case. For a general fibre $F$ of $X/S$ and $F''$ of $X''/S''$, we have by adjunction (cf. \thref{adjunction to fibre}):
   $$
   \mathrm{neft}(F;\mathcal{H}|_F) = \mathrm{psefft}(F;\mathcal{H}|_F) = \frac{1}{\mu} < \frac{1}{\mu ''} = \mathrm{psefft}(F'';\mathcal{H}''|_{F''}) = \mathrm{neft}(F'',\mathcal{H}''|_{F''})
   $$

     where $F$ and $F''$ are klt Fano varieties by adjunction. We have $-K_{F} \sim_{\mathbb{Q}} \frac{1}{\mu}\mathcal{H}|_{F}$, where $\mathcal{H}|_{F}$ is an ample integral linear system on $F$. Hence we can write $\frac{1}{\mu} = \frac{\mathrm{FI}(F)}{n}$, where $n$ is a positive integer and $\mathrm{FI}(F)$ is the Fano index of $F$. Hence by \thref{ACC for Fano index} in dimension $\leq d-1$, we have $\frac{1}{\mu}$ satisfies ACC and hence we cannot have infinitely many links of this type. Alternatively, one can use ACC for $\mathbb{R}$-complementary threshold (cf. \cite{Shokurov20}, Theorem 21).

   \item Suppose next we are in case 2-(IIIb) or 2-(IVb). In this case $\frac{1}{\mu} = \frac{1}{\mu ''}$ and the Sarkisov degree may not decrease. However, since we are running $(K_X + \frac{1}{\mu}\mathcal{H})$-MMP, and $(X,\frac{1}{\mu}\mathcal{H})$ is a canonical pair with $\mathcal{H}$ mobile, we conclude that the new Mori fibre space $X''/S''$ is still in Case 2, i.e. $c'' \geq \frac{1}{\mu ''}$. Hence it suffices to show that after finitely many links of type 2-(IIIb) and 2-(IVb), we must reach a link of type 2-(IVa) or type 2-(IIIa) and hence the pseudo-effective threshold strictly increase.

   \item \textbf{Case 2-(IIIb)}. For a link of type 2-(IIIb), we have $\mathrm{rank} \mathrm{Cl}(S'') = \mathrm{rank} \mathrm{Cl}(S) - 1$. Hence for a sequence of Sarkisov links of type 2-(IIIb) and 2-(IVb), the number of Sarkisov links of type 2-(IIIb) is bounded by $\mathrm{rank} \mathrm{Cl}(S)$.

   \item \textbf{Case 2-(IVb)}. Hence the remaining situation in case 2 is case 2-(IVb) (notice that this case would not happen for dimension $\leq 3$). Since the link is of type IV, every step appearing in the $(K_X + \frac{1}{\mu}\mathcal{H})$-MMP is a log flip. We claim further that the induced morphism on base
       $$
       \begin{tikzcd}
         S \arrow{rd} & & S'' \arrow{ld}\\
         & T &
       \end{tikzcd}
       $$
       is a $(K_S + \Delta_S)$-flip, where $(S,\Delta_S)$ is a choice of the adjunction pair as in \thref{Shokurov epsilon delta conjecture}.

       Indeed, by construction we have $\rho(S/T)=\rho(S''/T)=1$, $(K_S + \Delta_S)$ is anti-ample over $T$, and $(K_{S''} + \Delta_{S''})$ is ample over $T$. We first show that the diagram is small, i.e. isomorphic in codimension 1. Take a general curve $C$ in a fibre of $X/S$ such that $C$ doesn't intersect all the flipping locus of the sequence of log flips on the top. Let $C''$ be the strict transformation of $C$ in $X''$. Then
       $$
       0 = (K_X + \frac{1}{\mu}\mathcal{H})\cdot C = (K_{X''} + \frac{1}{\mu}\mathcal{H}'')\cdot C'' = 0
       $$
       Hence the birational map $X' \dashrightarrow X''$ induces a rational map $\phi : S \dashrightarrow S''$. The rational map $\phi$ is birational since the argument above can be inverted. This birational map must be small. Indeed, let $E$ be a prime divisor on $S$ contracted in $S''$. Then the divisor $f^{-1}(E)_{X''}$ is $\mathbb{Q}$-Cartier on $X''$. Let $C$ and $C''$ be the curve as above then we have $f^{-1}(E)_{X''} \cdot C'' =0$. Hence $f^{-1}(E)_{X''}$ is trivial over $S''$, but its center in $S''$ has codimension $\geq 2$ and $X''/S''$ is a Mori fibre space. Hence $f^{-1}(E)_{X''} \sim_{\mathbb{Q}} 0$, a contradiction.

       The morphism $S \rightarrow T$ is also birational. Otherwise, we can take a general curve $C$ in a fibre of $S/T$ such that it avoids the exceptional locus which has codimension $\geq 2$. Again we have
       $$
       0> (K_S + \Delta_S) \cdot C = (K_{S''} + \Delta_{S''}) \cdot C'' > 0
       $$
       which is a contradiction.

       Next we claim that we can choose $\Delta_{S''} = \phi_{*}\Delta_S$, that is, the log flip can be compatible with the adjunction. Indeed, we can choose a general ample $\mathbb{Q}$-divisor $A$ on $S$ such that $(K_X + \frac{1}{\mu}\mathcal{H} + f^{*}A)$ is trivial over $T$. Then the divisor $(K_{S} + \Delta_{S} + A)$ is also trivial over $T$. Hence the pair $(K_{X''} + \frac{1}{\mu}\mathcal{H}'' + \Phi_{*}f^{*}A)$ and $(K_{S''} + \phi_{*}\Delta_S + \phi_{*}A)$ is also trivial over $T$. Hence we have
       $$
       K_{X''} + \frac{1}{\mu}\mathcal{H}'' + \Phi_{*}f^{*}A \sim_{\mathbb{Q}} f^{\prime *}(K_{S''} + \phi_{*}\Delta_S + \phi_{*}A).
       $$ 
       Since $X/S$ and $X''/S''$ are Mori fibre spaces, the vertical divisor $\Phi_{*}f^{*}A$ is the pullback of a $\mathbb{Q}$-Cartier $\mathbb{Q}$-divisor on $S''$, which can only be $\phi_{*}A$. Hence we have $\Phi_{*}(f^{*}A) = f^{\prime *}\phi_{*}A$. We can then remove $A$ from both sides and conclude that $\Delta_{S''} = \phi_{*}\Delta_S$ satisfies the adjunction formula on $X''/S''$. Hence we conclude that the diagram on base is a $(K_S + \Delta_S)$-flip.

       By taking a general element $D$ in the $\mathbb{R}$-linear system $\frac{1}{\mu}\mathcal{H}$, we can assume that the coefficients of $D$ are $<\epsilon$ and $1-\epsilon <a(X,D) \leq 1$. Then the pair $(X,D)$ is $(1-\epsilon)$-lc. Hence the adjunction pair $(S,\Delta_S)$ can be chosen to have $\delta$-lc singularities by \thref{Shokurov epsilon delta conjecture}, where $\delta=\delta(d-d',1-\epsilon)$. Moreover, $\Delta_S$ is big by construction of the adjunction. The base locus of $\Delta_S$ lies inside the image of flipping locus. Hence by \thref{termination of delta lc flips} in dimension $d' \leq d-1$, such sequence must terminate.

       Alternatively, $(X,\frac{1}{\mu}\mathcal{H})$ has canonical singularities. By \thref{termination of canonical flips} in dimension $d$, the sequence of $(K_{X}+\frac{1}{\mu}\mathcal{H})$-flips on top must terminate.

   \item \textbf{Case 1.} Now suppose we are in Case 1 of \thref{construction of strong Sarkisov link}. Denote $p: \widetilde{X} \rightarrow X$ the maximal extraction of $X$ of an exceptional divisor $E$. Then by the condition of Case 1, $-(K_{\widetilde{X}} + c\mathcal{H}_{\widetilde{X}})$ is big over $S$. Hence we will always get a log Mori fibre space after running the $(K_{\widetilde{X}} + c\mathcal{H}_{\widetilde{X}})$-MMP. We claim that after running the log MMP, the pseudo-effective threshold $\frac{1}{\mu}$ is nondecreasing (One can further show that the equality holds if and only if either the dimension of the base $S$ increase, or the link is square, cf. \thref{square map}). Indeed, by our assumption $c < \frac{1}{\mu}$, i.e. the pair $(X,\frac{1}{\mu}\mathcal{H})$ is not canonical. Hence we have $p^{*}(K_{X} + \frac{1}{\mu}\mathcal{H}) = K_{\widetilde{X}} + \frac{1}{\mu}\widetilde{\mathcal{H}} + aE$ for some positive real number $a$. Take a general curve $C''$ in the general fibre of $X'' \rightarrow S''$, such that $C''$ doesn't intersect all the exceptional locus and flipped locus of the link. Then we have:
   \begin{align*}
   (K_{X''}+\frac{1}{\mu}\mathcal{H}'')\cdot C'' &= (K_{\widetilde{X}}+\frac{1}{\mu}\widetilde{\mathcal{H}})\cdot \tilde{C}\\
   &\leq (K_{\widetilde{X}} + \frac{1}{\mu}\widetilde{\mathcal{H}} + aE)\cdot \widetilde{C}\\
   &=(K_{\tilde{X}} + \frac{1}{\mu}\mathcal{H})\cdot C = 0
   \end{align*}
   Hence we have $\frac{1}{\mu} \leq \frac{1}{\mu''}$.

   Next we look into the canonical threshold. Since we are running the $(K_{\widetilde{X}}+c\widetilde{\mathcal{H}})$-MMP and $\widetilde{\mathcal{H}}$ is mobile, the pair $(X'',c\mathcal{H}'')$ is canonical by \thref{prop:singularity after MMP}. Hence the canonical threshold cannot decrease. Now suppose the canonical thresholds are equal. We want to show that the number of crepant divisors strictly decrease, and hence the Sarkisov degree always strictly decrease. Indeed, assume the link is of type I, then the exceptional divisor $E$ of the maximal extraction becomes a divisor on $X''$, so the number of crepant divisors decrease by at least 1 in this case. Next assume that the link is of type II. Since the divisorial contraction to $X''$ is a $(K+c\mathcal{H})$-negative extremal contraction, the exceptional divisor $E''$ is a crepant divisor of $(X,c\mathcal{H})$, but cannot be a crepant divisor of $(X'',c\mathcal{H}'')$. Hence the number of crepant divisors also decrease by at least 1 in this case.

   In every situation of Case 1, the Sarkisov degree strictly decrease. By the discussion in Case 2-(IIIa) and 2-(IVa), the pseudo-effective thresholds satisfy ACC. It is obvious that the number of crepant divisors satisfies DCC. However, ACC for canonical thresholds remains open. To finish the proof, we need to show that there is no infinite sequence of Sarkisov links such that the pseudo-effective threshold is stationary and the canonical threshold strictly increase. Suppose to the contrary that there exists an infinite sequence
   $$
   \begin{tikzcd}
     X_0 \arrow[dashed,r] \arrow{d} & X_1 \arrow[dashed,r] \arrow{d} & \cdots\\
     S_0 \arrow[dashed,r]  & S_1 \arrow[dashed,r]  & \cdots
   \end{tikzcd}
   $$
   of Sarkisov links of Case 1, such that the pseudo-effective threshold $\frac{1}{\mu}$ is stationary but the canonical threshold is not. Since $c_i < \frac{1}{\mu}$ is bounded, the canonical thresholds $c_i$ increase to a limit $\alpha \leq \frac{1}{\mu}$.

   We proceed in several steps to show this is impossible.
   \begin{enumerate}[label = Step \arabic* :,align=left]
     \item Let $(Z_i, \alpha \mathcal{H}_{Z_i})$ be the maximal extraction of the divisor $E_i$ corresponding to the link $X_i \dashrightarrow X_{i+1}$. We show that the pair $(X_i, \alpha \mathcal{H}_{X_i})$ and $(Z_i, \alpha \mathcal{H}_{Z_i})$ have only log canonical singularities for $i \gg 0$.

         Indeed, suppose to the contrary that $c_i \leq \mathrm{lct}(X_i, \mathcal{H}_{X_i}) < \alpha$ for infinitely many $i$. Then we have a subsequence such that $\mathrm{lct}(X_i, \mathcal{H}_{X_i})$ strictly increase but not stabilize. This contradicts \thref{ACC for log canonical thresholds}. The same argument works for $(Z_i, \alpha \mathcal{H}_{Z_i})$.
     \item We show that the $(K_{Z_{i}} + c_i \mathcal{H}_{Z_{i}})$-flips on $Z_i$ are also $(K_{Z_{i}} + \alpha \mathcal{H}_{Z_{i}})$-flips. As a consequence, for any divisor $E$ over $X$, we have
         $$
         a(E, X_0, \alpha \mathcal{H}_{X_0}) \leq a(E, X_i, \alpha \mathcal{H}_{X_i})
         $$
         and the equality holds if and only if the composition of links $X_0 \dashrightarrow X_i$ is an isomorphism at the center of $E$.

         Indeed, this sequence of flips is obtained by running the 2-ray game on $Z_i/S_i$. Notice that $K_{Z_{i}} + \alpha \mathcal{H}_{Z_{i}} \equiv_{S_i} -a_i E_i$ for some $0 < a_i \leq 1$ since $\alpha > c_i$ and the pair $(X_i, \alpha \mathcal{H}_{X_i})$ is lc. We conclude that every $(K_{Z_{i}} + c_i \mathcal{H}_{Z_{i}})$-flipping contraction must also be $(K_{Z_{i}} + \alpha \mathcal{H}_{Z_{i}})$-negative, otherwise there exists a log minimal model of $(Z_{i},\alpha \mathcal{H}_{Z_{i}})$ over $S_i$, but $(K_{Z_{i}} + \alpha \mathcal{H}_{Z_{i}})$ is not pseudo-effective over $S_i$, a contradiction. The consequence comes from the fact that
         $$
         a(E,X_{i},\alpha \mathcal{H}_{X_{i}}) = a(E,Z_{i},\alpha \mathcal{H}_{Z_{i}} + a_{i}E_{i}) \leq a(E,Z_{i},\alpha \mathcal{H}_{Z_{i}}) \leq a(E, X_{i+1}, \alpha \mathcal{H}_{X_{i+1}})
         $$
         for every link in the sequence.

     \item The pair $(X_i, \alpha \mathcal{H}_{X_i})$ is plt for $i \gg 0$.

     Indeed, suppose to the contrary that there exist infinitely many pairs $(X_i, \alpha \mathcal{H}_{X_i})$ with lc but not plt singularities. Then there exists exceptional divisors $E_i$ over $X_{i}$ and the log discrepancy
     $$
     a(E_i, X_0, \alpha \mathcal{H}_{X_0}) = a(E_i, X_i, \alpha \mathcal{H}_{X_i}) = 0.
     $$
     Hence by the previous step, the composition of links $X_0 \dashrightarrow X_i$ is an isomorphism at the center of $E_i$. Then the local canonical thresholds
     $$
     c(\eta_{\mathrm{Center}_{X_i}(E_i)}\in X_i,\mathcal{H}_{X_i}) = c(\eta_{\mathrm{Center}_{X_0}(E_i)}\in X_0,\mathcal{H}_{X_0}).
     $$ We also have the inequality
     $$
     c_i \leq c(\eta_{\mathrm{Center}_{X_i}(E_i)} \in X_i,\mathcal{H}_{X_i})  < \alpha.
     $$
     But on the other hand, the set $\{ c(x \in X_0,\mathcal{H}_{X_0})|x \in X_0 \}$ is finite by \thref{lower semi-continuity of canonical threshold}, a contradiction.

     \item We assume without loss of generality that all pairs $(X_i, \alpha \mathcal{H}_{X_i})$ are plt. For every exceptional divisor $E_i$ of a maximal extraction, we have
         $$
         a(E_i,X_0,\alpha \mathcal{H}_X) \leq a(E_i,X_i,\alpha \mathcal{H}_{X_i}) < 1.
         $$
         All $E_i$ are distinct prime divisors. Hence we have infinitely many divisors $E$ over $X_0$ such that $a(E, X_0,\alpha\mathcal{H}_{X_0}) < 1$, contradicting \thref{finiteness of divisors with negative discrepancy}.
   \end{enumerate}

   \end{enumerate}
\end{proof}

\begin{corollary}\label{cor:strong Sarkisov program in dim 3 and 4}
  The strong Sarkisov program holds in dimension 3 and 4.
\end{corollary}

\begin{proof}
  There is no flip in dimension 2. Termination of 4-fold canonical flips is proved in \cite{Fuj03}.
\end{proof}

\begin{corollary}
  Assume termination of $\frac{1}{2}$-lc flips with big boundary in dimension 4. Then the strong Sarkisov program holds in dimension 5.
\end{corollary}

\begin{proof}
  It was shown in \cite{HJL}, Corollary 1.8 that we can choose boundary $\Delta_{S}$ such that $(X,\Delta_{S})$ is $\frac{1}{2}$-lc.
\end{proof}

  The next proposition shows that at least a very weak version of termination of flips is needed in order to obtain termination of strong Sarkisov program.

\begin{proposition}\label{prop:strong Sarkisov implies termination of flips}
  Assume that the strong Sarkisov program in dimension $d$ holds. Then there is no infinite sequence of $(K_{X}+\Delta)$-flips
  $$
  (X,\Delta) \dashrightarrow (X_{1},\Delta_{X_{1}}) \dashrightarrow \cdots
  $$
  where $(X,\Delta)$ is canonical, $X$ is terminal of dimension $\leq d-1$ and $(K_X + \Delta)$ big and mobile.
\end{proposition}

\begin{proof}
  Suppose to the contrary that there exists a infinite sequence of $(K_X + \Delta_X)$-flips in dimension $d-1$ satisfying the conditions. Then we can construct an infinite sequence of Sarkisov links of type $\mathrm{IV}_s$ in dimension $d$:
  $$
  \begin{tikzcd}
    X_0 \times \mathbb{P}^1 \arrow[dashed,rr]\arrow{d} && X_1 \times \mathbb{P}^1 \arrow[dashed,rr]\arrow{d} && \dots \\
    X = X_0 \arrow[dashed,rr] \arrow{rd} && X_1 \arrow[dashed,rr] \arrow{ld} \arrow{rd} && \dots \arrow{ld} \\
    & Y_0 && Y_1 &
  \end{tikzcd}
  $$
    By our assumption $(K_X + \Delta_X)$ is big and mobile. Since $(K_{X}+\Delta)$ is big, the log canonical model of $(X,\Delta)$ exists. Let $(X',\Delta_{X'})$ be the log canonical model of $(X,\Delta)$, and fix the Mori fibre space $X' \times \mathbb{P}^1 / X'$. Take $\mathcal{H}' = |m( 2 X' \times \{ pt \} + (\Delta_{X'}) \times \mathbb{P}^1 + n (K_{X'} + \Delta_{X'}) \times \mathbb{P}^1)|$ for $m \gg 0$ and $(n+1)$ sufficiently divisible. We have $\frac{1}{\mu} =\frac{1}{\mu '}= \frac{1}{m}$ is stable, and $K_{X' \times \mathbb{P}^1} + \frac{1}{m}\mathcal{H}' \sim (n+1)(K_{X'}+ \Delta_{X'})$. Hence the condition in \thref{Sarkisov degree} is satisfied. However, the strong Sarkisov program of the birational map $X \times \mathbb{P}^1 \dashrightarrow X' \times \mathbb{P}^1$ does not terminate in the above diagram, a contradiction.
\end{proof}

  As a byproduct of the proof of \thref{termination of strong Sarkisov link}, we can give another proof of the weak Sarkisov program:

  \begin{theorem}\label{theorem: weak Sarkisov program with decreasing Sarkisov degree}
    Notation as in \thref{Sarkisov degree}. Then $\Phi$ can be factorized into a finite sequence of elementary links by the strong Sarkisov program.
  \end{theorem}

  \begin{proof}
    The proof is the same as the proof of \thref{termination of strong Sarkisov link}, except that we make a special choice of the flip in the case 2-(IVb).

    Assume that we have the birational map
  $$
  \begin{tikzcd}
    X \arrow[dashed,r,"\Phi"]\arrow[d,"f"] & X' \arrow[d,"f'"] \\
    S & S'
  \end{tikzcd}
  $$
  where $\frac{1}{\mu} \leq c$. Take a general very ample divisor $A$ on $S$. We run the $(K_S + \Delta_S)$-MMP on $S$ \emph{with scaling of $A$} to obtain an extremal contraction $S \rightarrow Y$. We claim that this construction coincide with the construction in Case 2 of \thref{construction of strong Sarkisov link}. Indeed, by definition of MMP with scaling, the extremal ray $R'_{S}$ on $S$ is chosen so that the number $a_{S} = \frac{-(K_{S}+\Delta_{S})\cdot v'_{S}}{A \cdot v'_{S}}>0$ is the largest. Let $f^{*}A$ be the supporting hyperplane of the extremal ray $R$ on $X$, and $R'$ be the extremal ray such that the number $a = \frac{-(K_X + \frac{1}{\mu}\mathcal{H}) \cdot v'}{(f^{*}A) \cdot v'} > 0$ is the largest. Hence by the projection formula, we can take $R'_{S} = f_{*}R'$.

   If the resulting link is not of type 2-(IVb), then by the proof of \thref{termination of strong Sarkisov link} the sequence terminates. Otherwise we have a link of type 2-(IVb):
    $$
  \begin{tikzcd}
    X = X_0 \arrow[dashed,rr,"\Phi"]  \arrow{d} && X_1 \arrow[d] \\
    S = S_0 \arrow[dashed,rr,"\phi"] \arrow{rd} && S_1 \arrow[ld] \\
    & Y = Y_0 &
  \end{tikzcd}
  $$
    Now $\phi_{*}A$ is not ample and hence doesn't define a supporting hyperplane. However, we claim that the same construction still gives us an extremal ray $R'_{1}$ such that if we denote by $R_{1}$ the extremal ray of $X_{1}/S_{1}$, then $R_{1},R'_{1}$ spanned a 2 dimensional extremal face. Indeed, let $R'_{1}$ be the extremal ray on $X_{1}$ such that the number $a_{1} = \frac{-(K_{X_{1}} + \frac{1}{\mu}\mathcal{H}_{1}) \cdot v'_{1}}{(f^{*}(\phi_{*}A)) \cdot v'_{1}} > 0$ is the largest. Such extremal ray is unique for sufficiently general $A$. Consider the hyperplane $(K_{X_{1}} + \frac{1}{\mu}\mathcal{H}_{1} + a_{1}(f^{*}(\phi_{*}A))\cdot v_{1}=0$. This hyperplane is a supporting hyperplane of the face spanned by $R_{1}$ and $R'_{1}$. Indeed, if the ray is not $(K_{X_{1}} + \frac{1}{\mu}\mathcal{H}_{1})$-negative and not in the flipped locus, then the intersection is always positive. Now for any other $(K_{X_{1}} + \frac{1}{\mu}\mathcal{H}_{1})$-negative extremal ray the intersection is positive by construction. Since $A$ is taken to be sufficiently general, we have $a_{1} < a$. Hence for any flipped locus the intersection is also positive.

    Hence we can continue the MMP with scaling of $A$ on the base, and by \cite{BCHM}, the sequence terminates and we will not have infinite sequence of links of type 2-(IVb).
  \end{proof}

  \section{Decreasing Sarkisov degree}

  Another problem in the strong Sarkisov program is the Sarkisov degree may not decrease as we originally expected. We want to introduce an augmented version of the Sarkisov degree, so that it decreases in every step. Firstly we recall the definition of difficulties.
  
  \begin{definition}[Difficulties, cf. \cite{AHK}]\label{def: difficulties}
  \hspace{2em}
    \begin{itemize}
      \item \textbf{Weight functions.} We take a weight function $w: (0,+\infty) \rightarrow \mathbb{R}_{\geq 0}$ satisfying the following conditions:
          \begin{enumerate}
            \item $w(a) \geq 0$ for all $a$, and $w(a) = 0$ for $a \geq 1$.
            \item $w$ is decreasing (not necessarily strictly).
            \item Whenever $0 < \sum m_{i}b_{i} < 1$ for some $m_{i} \in \mathbb{Z}_{\geq 0}$, one has
            $$
            w(1 - \sum m_{i}b_{i}) \geq \sum m_{i}w(1-b_{i}).
            $$
            \item For $i=1,\cdots,m$, let $\{b_{i}^{n}, n \in \mathbb{N} \}$ be $m$ decreasing sequences of real numbers in $(0,1)$. Then there exists $n_{0}$ such that the following holds: Whenever $k \geq 2, m_{i} \geq 0$ and $n \geq n_{0}$ are integers such that
                $$
                1 - \sum m_{i}b_{i}^{n} \in (0,1), \quad \mathrm{resp.} \quad k(1 - b_{i}^{n}) - \sum\limits_{i \neq i_{0}} m_{i}b_{i}^{n} \in (0,1),
                $$
                then one has
                $$
                \mathrm{diff}_{1}(n) \geq \mathrm{diff}_{1}(n+1), \quad \mathrm{resp.} \quad \mathrm{diff}_{2}(n) \geq \mathrm{diff}_{2}(n+1),
                $$
                where
                $$
                \mathrm{diff}_{1}(n) = w(1 - \sum m_{i}b_{i}^{n}) - \sum m_{i}w(1-b_{i}^{n})
                $$
                and
                $$
                \mathrm{diff}_{2}(n) = w (k(1 - b_{i}^{n}) - \sum\limits_{i \neq i_{0}} m_{i}b_{i}^{n}) - w(k(1-b_{i}^{n})).
                $$
                For example, we could take for $w$ a piecewise linear function with decreasing absolute values of slopes. Let $\alpha \in (0,1)$. We define
                $$
                w_{\alpha}^{-}(x) = \begin{cases}
                                   1-x, & \mbox{if } x \leq \alpha \\
                                   0, & \mbox{otherwise}
                                 \end{cases}
                $$
                and
                $$
                w_{\alpha}^{+}(x) = \begin{cases}
                                   1-x, & \mbox{if } x < \alpha \\
                                   0, & \mbox{otherwise}.
                                 \end{cases}
                $$
          \end{enumerate}
      \item \textbf{Summed weight functions.} Let $w$ be a weight function. We define the \emph{summed weight function} $W: (-\infty,1) \rightarrow \mathbb{R}_{\geq 0}$ by the formula
          $$
          W(b) = \sum\limits_{k=1}^{+\infty}w(k(1-b)).
          $$
          The sum makes sense since there are only finitely many nonzero terms for every $b \in (-\infty,1)$.
      \item \textbf{Difficulties, terminal case.} Let $(Y,B = \sum b_{i}B_{i})$ be a terminal pair, not necessary effective, i.e. possibly with some coefficients $b_{i}<0$. Let $\nu: \coprod \widetilde{B_{i}} \rightarrow \bigcup B_{i}$ be the normalization of the divisor $\mathrm{Supp}B$. For any irreducible subvariety $C \subseteq \mathrm{Supp}B$ the preimage $\widetilde{C} = \nu^{-1}(C)$ splits into a union of irreducible components $\widetilde{C}_{ij} \subseteq \widetilde{B}_{i}$.

          Then we define \emph{the difficulty of $(Y,B)$} to be
          $$
          d(Y,B) = \sum\limits_{a(B_{j})\leq 0} W(b_{i})\rho(\widetilde{B}_{i}) + \sum\limits_{v;\mathrm{codim}C(v,Y) \neq 2} w(a_{v}) + \sum\limits_{\mathrm{irr}.C \subseteq Y; \mathrm{codim}C=2}\left[\sum\limits_{v;C(v,Y)=C}w(a_{v}) + \sum\limits_{\widetilde{C}_{ij}}W(b_{i})\right].
          $$
      \item \textbf{Difficulties, klt case.} Let $(X,D)$ be a klt pair. The difficulty of $(X,D)$ is defined to be $d(X,D) = d(Y,B)$, where $f: (Y,B) \rightarrow (X,D)$ is any terminal pair such that $K_{Y}+B = f^{*}(K_{X}+D)$.
    \end{itemize}
  \end{definition}

  \begin{definition}[Augmented Sarkisov degree]\label{def:augmented Sarkisov degree}
    Notation as in \thref{Sarkisov degree}. We define the \emph{augmented Sarkisov degree} to be the 6-tuple $(\mu,b,\rho,d,c',e')$, where:
    \begin{enumerate}[align=left]
      \item If $\frac{1}{\mu} \leq c$ then we set $c' = +\infty$ and $e'=0$, and otherwise $c'=c$ and $e'=e$;
      \item $\rho = \mathrm{rank} \mathrm{Cl}(S)$ is the rank of the class group of the base $S$, we set $\rho = +\infty$ if $\frac{1}{\mu} > c$;
      \item $d = d(S,\Delta_S)$ is the difficulty of the pair $(S,\Delta_S)$. We let $d = +\infty$ if $b \leq 2$ or $\frac{1}{\mu} > c$.

      The order ``$<$'' is defined to be the lexicographic order of $(<,>,<,<,>,<)$.
    \end{enumerate}
  \end{definition}

  \begin{proposition}\label{prop:decrease Sarkisov degree}
    The augmented Sarkisov degree decreases in every step of a Strong Sarkisov program, and strictly decreases except for Case 2-(IVb).
  \end{proposition}

  \begin{proof}
    One can just verify by looking at every situation in the proof of \thref{termination of strong Sarkisov link}. More precisely, we can trace the change of the augmented Sarkisov degree by the following table:
    \begin{center}
    \begin{longtable}{|C|C|C|C|C|C|}
    \hline
    \mathrm{Case} &\frac{1}{\mu}  & \rho & d & c' & e' \\
    \hline
    1 & \nearrow ? & +\infty \searrow? & +\infty \searrow? & \nearrow ? & \searrow\\
    \hline
    \mathrm{2-(IIIa)} & \nearrow  & / & / & / & / \\
    \hline
    \mathrm{2-(IIIb)} & = & \searrow & / & / & / \\
    \hline
    \mathrm{2-(IVa)} & \nearrow & / & / & / & / \\
    \hline
    \mathrm{2-(IVb)} & = & = & \searrow ? & +\infty & 0\\
    \hline
    \end{longtable}
    \end{center}
    Here a question mark after an arrow means that the change is not strict, and we look at the arrow in the next column only when the equality holds.
  \end{proof}
  
  \begin{statement}[strictness of decreasing]
  The augmented Sarkisov degree may not strictly decrease exactly because of the fact that the difficulty may not strictly decrease after a log flip. One can fix this by finding some invariants which strictly decrease for every log flip.
  \end{statement}
  
  \begin{corollary}\label{cor:factorization with decreasing degree}
    Let $\Phi: X/S \dashrightarrow Y/T$ be a birational map between Mori fibre spaces. Then there exists a factorization of $\Phi$ into Sarkisov links, such that the augmented Sarkisov degree decreases in every step.
  \end{corollary}
  
  \begin{proof}
    Immediate from \thref{theorem: weak Sarkisov program with decreasing Sarkisov degree} and \thref{prop:decrease Sarkisov degree}.
  \end{proof}

  \section{Relation to the syzygy of Mori fibre spaces}

  \begin{convention}
    In this section, we assume that the strong Sarkisov program holds.
  \end{convention}

  We recall the homotopical syzygies of Mori fibre spaces (cf. \cite{Myself_Syzygy}):

    \begin{theorem}[Homotopical syzygies of $Y/R$]\label{homotopical syzygy}
    There exists a regular (possibly infinite dimensional) CW complex $CW = CW(Y/R)$ satisfying the following:
  \begin{enumerate}[label=(\arabic*),align=left]
    \item If $CW$ is not the empty set, then $CW$ is contractible.
    \item There is a one-to-one correspondence between cells of $CW$ of dimension $d$, and central models of rank $r = d + 1$ of $Y/R$.
    \item For any 2 cells $C_{1}$ and $C_{2}$, the inclusion $C_{1} \subseteq \overline{C_{2}}$ holds if and only if the corresponding rank $r_{i}$ central models satisfy $X_{1}/Z_{1} \preceq X_{2}/Z_{2}$.
    \item There exists a natural action of $G$ on the set of all $d$-cells of $CW$ for any $d \geq 0$. For any $d$-cell $C_{d}$ corresponding to a central model $X_{r}/Z_{r}$ of rank $r$, the stabilizer of $\left[ C_{d} \right]$ is $\mathrm{Aut}(X_{r}\rightarrow Z_{r}/R)$, the fibrewise regular automorphism group of $X_{r} \rightarrow Z_{r}$ over $R$.
  \end{enumerate}
  \end{theorem}

  In this article we set $R$ to be a point. Then the strong Sarkisov program implies:

  \begin{theorem}\label{thm:Sarkisov diagram}
    Let $X'/S'$ be a Mori model of $Y$. Fix a very ample complete linear system $\mathcal{H}' = |m(-K_{X'} + f^{\prime *}A')|$ on $X'$, where $m$ is a sufficiently large integer and $A'$ is a very ample divisor on $S'$. Then for any Mori model $X/S$ of $Y$, there exists a directed graph $\mathfrak{G}(X/S,X'/S',\mathcal{H}')$ such that:
    \begin{enumerate}[label=(\arabic*)]
      \item The undirected graph of $\mathfrak{G}(X/S,X'/S',\mathcal{H}')$ is a finite subgraph of $CW_{1}(Y)$.
      \item The vertex corresponding to $X/S$ is the unique source, and the vertex corresponding to $X'/S'$ is the unique destination of $\mathfrak{G}(X/S,X'/S',\mathcal{H}')$.
      \item There is a 1-1 correspondence
      \begin{align*}
      \{\text{Paths of } \mathfrak{G} \text{ from the vertex corresponding to } X/S \text{ to the vertex corresponding to } X'/S'\} \\
       \longleftrightarrow \{ \text{Sequences of untwisting in the strong Sarkisov program of } X/S \dashrightarrow X'/S'\}
      \end{align*}
      and an injection
      \begin{align*}
      \{ \text{Steps of intermediate untwisting in the strong Sarkisov program of } X/S \dashrightarrow X'/S' \} \\
       \hookrightarrow \{\text{Edges of }\mathfrak{G}(X/S,X'/S',\mathcal{H}') \}.
      \end{align*}
    \end{enumerate}
  \end{theorem}

  \begin{proof}
    We construct the graph from \thref{construction of strong Sarkisov link}. We start from the vertex corresponding to $X/S$ in $CW(Y)$. For every possible untwisting $X/S \dashrightarrow X_{1}/S_{1}$ we draw an arrow (i.e. a directed edge) from the vertex corresponding to $X/S$ to the vertex corresponding to $X_{1}/S_{1}$, and we repeat the procedure for $X_{1}/S_{1}$.

    It remains to show the finiteness of the above construction. Firstly, for any vertex $V$, the number of arrows sourced from $V$ is finite. Indeed, in Case 1 of \thref{construction of strong Sarkisov link} the number of crepant divisors of the pair $(X,c\mathcal{H})$ is finite, and in Case 2 the number of negative adjacent extremal rays is also finite. Secondly, the lengths of paths are bounded. Indeed, suppose the lengths are not bounded, then among all the possible untwisting $X/S \dashrightarrow X_{1}/S_{1}$, the exists a Mori model $X_{1}/S_{1}$ such that the lengths of paths from $X_{1}/S_{1}$ to $X'/S'$ is unbounded. Repeat this construction we can find an infinite sequence of untwisting, contradicting the strong Sarkisov program.
  \end{proof}

\printbibliography

\end {spacing}
\end {document}